\documentclass[11pt,reqno]{amsart}
\usepackage{amsfonts}
\usepackage{amssymb}
\usepackage{amsmath}
\usepackage{eepic}
\usepackage{dsfont}

\allowdisplaybreaks

\newtheorem{theorem}{Theorem}
\newtheorem{lemma}[theorem]{Lemma}

\theoremstyle{definition}

\newtheorem{remark}[theorem]{Remark}
\theoremstyle{plain}

\newcommand{\Ac}{\mathcal A}

\newcommand{\GG}{\mathds G}

\newcommand{\Bc}{\mathcal B}

\newcommand{\inv}{^{-1}}

\def\arr#1{\mathop{{\buildrel #1\over\longrightarrow}}}
\def\ov{\overline}
\def\wh{\widehat}

\def\power#1{\langle\!{\langle}#1\rangle\!\rangle}
\def\inner#1#2{\langle #1,#2\rangle}

\begin{document}
    \title[Testing spherical transitivity in iterated wreath products]{Testing spherical transitivity in iterated wreath products of cyclic groups}
    \keywords{Automata, spherical transitivity, iterated wreath products, rooted trees, rational power series}

\author{Benjamin Steinberg}\address{School of Mathematics and
  Statistics\\ Carleton
University\\ 1125 Colonel By Drive\\ Ottawa\\ Ontario  K1S 5B6 \\
Canada}\email{bsteinbg@math.carleton.ca}
\thanks{The author acknowledges the support of NSERC}
\date{Version of \today}

\begin{abstract}
We give a partial solution a question of Grigorchuk, Nekrashevych,
Sushchanskii and \v Suni\'k by giving an algorithm to test whether a
finite state element of an infinite iterated (permutational) wreath
product $\wh G = \mathbb Z/k\mathbb Z\wr \mathbb Z/k\mathbb Z\wr
\mathbb Z/k\mathbb Z\wr \cdots$ of cyclic groups of order $n$ acts
spherically transitively. We can also decide whether two finite
state spherically transitive elements of $\wh G$ are conjugate. For
general infinite iterated wreath products, an algorithm is
presented to determine whether two finite state automorphisms have
the same image in the abelianization.
\end{abstract}
\maketitle 

\section{Introduction and main results}
The purpose of this note is to offer a partial solution to a
question of Grigorchuk, Nekrashevych, Sushchanskii and \v
Suni\'k~\cite{GNS,selfsimilarsurvey}. Let $T_k$ be the rooted
regular $k$-ary tree.  We view it as the Cayley graph of the free
monoid $A_k^*$, where $A_k=\{0,\ldots,k-1\}$ is the standard
alphabet of size $k$.  In particular, we identify vertices with
words.  It is well known that $\mathrm{Aut}(T_k)$ is a profinite
group. In fact, there is a permutational wreath product
decomposition $(\mathrm{Aut}(T_k),T_k) = (S_k,A_k)\wr
(\mathrm{Aut}(T_k),T_k)$~\cite{branch,Bass,GNS}. Thus
$\mathrm{Aut}(T_k) = (S_k,A_k)\wr (S_k,A_k)\wr \cdots$.  For more on
this group see~\cite{branch,Bass,GNS,selfsimilarsurvey,selfsimilar}.
An element $f\in \mathrm{Aut}(T_k)$ is said to be \emph{spherically
transitive} if, for each $n$,  $\langle f\rangle$ acts transitively
on the set of vertices at distance $n$ from the root, i.e.\
transitively on the set of words of length
$n$~\cite{branch,Bass,GNS,selfsimilarsurvey,selfsimilar}. This is
equivalent to topological transitivity and ergodicity of the action
on the boundary  $\partial T_k$~\cite{GNS}.

If $f\in \mathrm{Aut}(T_k)$ has wreath product decomposition
$\lambda_f(f_0,\ldots,f_{k-1})$, then $f_i$ is called the
\emph{section} of $f$ at $i\in A_k$.  We shall use the notation
$\lambda_f$ throughout for the element of $S_k$ associated to $f$.
One can the define inductively, for any word $w\in A_k^*$, the
section $f_w$ by the formula $f_{ua} = (f_u)_a$ where $a\in A_k$ and
$u\in A_k^*$. Of course, $f_{\varepsilon} =f$, where $\varepsilon$
is the empty word. One then has the formula $f(uw) = f(u)f_u(w)$ for
any words $u,w\in A_k^*$.  An element $f\in \mathrm{Aut}(T_k)$ is
said to be \emph{finite state} if it has only finitely many distinct
sections. This is the same as saying that $f$ can be computed by a
finite state automaton.

A \emph{finite state automaton} over an alphabet $A$ is a $4$-tuple $\Ac =
(Q,A,\delta,\lambda)$ where $Q$ is a finite set of states,
$\delta:Q\times A\to A$ is the transition function and
$\lambda:Q\times A\to A$ is the output function.  We set $q_a =
\delta(q,a)$ and $q(a) = \lambda(q,a)$ for $q\in Q$, $a\in A$.  We
extend this to words by the formulas:
\begin{gather} q_{au} = (q_a)_u
  \\\label{action} q(au) = q(a)q_a(u)
\end{gather}
So each state $q\in \mathcal A$ gives rise to a function, via
\eqref{action}, from $A^*\to A^*$ (in fact an endomorphism of the
rooted Cayley tree of $A^*$), also denoted by $q$.  An automaton
with a distinguished state is called an \emph{initial automaton}.

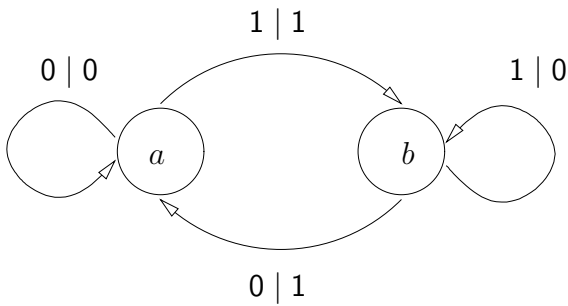
\begin{figure}[htbp]
\begin{center}
\setlength{\unitlength}{0.00083333in}
\begingroup\makeatletter\ifx\SetFigFont\undefined%
\gdef\SetFigFont#1#2#3#4#5{%
  \reset@font\fontsize{#1}{#2pt}%
  \fontfamily{#3}\fontseries{#4}\fontshape{#5}%
  \selectfont}%
\fi\endgroup%
\begin{picture}(3444,1800)(0,-10)
\put(1722.000,450.000){\arc{2121.320}{3.9270}{5.4978}}
\whiten\path(2362.476,1257.485)(2472.000,1200.000)(2402.413,1302.263)(2362.476,1257.485)
\put(1722.000,1350.000){\arc{2121.320}{0.7854}{2.3562}}
\whiten\path(1081.524,542.515)(972.000,600.000)(1041.587,497.737)(1081.524,542.515)
\put(972,900){\ellipse{540}{540}} \put(2472,900){\ellipse{540}{540}}
\path(2757,810)(2758,808)(2762,805)
        (2767,798)(2776,788)(2788,776)
        (2802,760)(2819,743)(2837,724)
        (2858,705)(2879,685)(2902,666)
        (2926,648)(2951,632)(2977,617)
        (3004,604)(3034,594)(3065,587)
        (3098,584)(3132,585)(3160,590)
        (3186,597)(3210,606)(3232,616)
        (3252,626)(3269,636)(3284,646)
        (3297,656)(3309,665)(3320,674)
        (3329,682)(3338,691)(3347,700)
        (3356,710)(3364,720)(3373,732)
        (3382,744)(3391,759)(3401,775)
        (3410,794)(3418,814)(3425,836)
        (3430,860)(3432,885)(3430,910)
        (3425,934)(3418,956)(3410,976)
        (3401,995)(3391,1011)(3382,1026)
        (3373,1038)(3364,1050)(3356,1060)
        (3347,1070)(3338,1079)(3329,1088)
        (3320,1096)(3309,1105)(3297,1114)
        (3284,1124)(3269,1134)(3252,1144)
        (3232,1154)(3210,1164)(3186,1173)
        (3160,1180)(3132,1185)(3098,1186)
        (3065,1183)(3034,1176)(3004,1166)
        (2977,1153)(2951,1138)(2926,1122)
        (2902,1104)(2879,1085)(2858,1065)
        (2837,1046)(2819,1027)(2802,1010)
        (2788,994)(2776,982)(2757,960)
\whiten\path(2812.730,1070.427)(2757.000,960.000)(2858.139,1031.210)(2812.730,1070.427)
\path(687,990)(686,992)(682,995)
        (677,1002)(668,1012)(656,1024)
        (642,1040)(625,1057)(607,1076)
        (586,1095)(565,1115)(542,1134)
        (518,1152)(493,1168)(467,1183)
        (440,1196)(410,1206)(379,1213)
        (346,1216)(312,1215)(284,1210)
        (258,1203)(234,1194)(212,1184)
        (192,1174)(175,1164)(160,1154)
        (147,1144)(135,1135)(124,1126)
        (115,1118)(106,1109)(97,1100)
        (88,1090)(80,1080)(71,1068)
        (62,1056)(53,1041)(43,1025)
        (34,1006)(26,986)(19,964)
        (14,940)(12,915)(14,890)
        (19,866)(26,844)(34,824)
        (43,805)(53,789)(62,774)
        (71,762)(80,750)(88,740)
        (97,730)(106,721)(115,712)
        (124,704)(135,695)(147,686)
        (160,676)(175,666)(192,656)
        (212,646)(234,636)(258,627)
        (284,620)(312,615)(346,614)
        (379,617)(410,624)(440,634)
        (467,647)(493,662)(518,678)
        (542,696)(565,715)(586,735)
        (607,754)(625,773)(642,790)
        (656,806)(668,818)(687,840)
\whiten\path(631.270,729.573)(687.000,840.000)(585.861,768.790)(631.270,729.573)
\put(222,1350){\makebox(0,0)[lb]{\smash{{{\SetFigFont{12}{14.4}{\rmdefault}{\mddefault}{\updefault}$\mathsf{0}\mid
\mathsf{0}$}}}}}
\put(3147,1350){\makebox(0,0)[lb]{\smash{{{\SetFigFont{12}{14.4}{\rmdefault}{\mddefault}{\updefault}$\mathsf{1}\mid
\mathsf{0}$}}}}}
\put(1522,1650){\makebox(0,0)[lb]{\smash{{{\SetFigFont{12}{14.4}{\rmdefault}{\mddefault}{\updefault}$\mathsf{1}\mid
\mathsf{1}$}}}}}
\put(1522,0){\makebox(0,0)[lb]{\smash{{{\SetFigFont{12}{14.4}{\rmdefault}{\mddefault}{\updefault}$\mathsf{0}\mid
\mathsf{1}$}}}}}
\put(897,825){\makebox(0,0)[lb]{\smash{{{\SetFigFont{12}{14.4}{\rmdefault}{\mddefault}{\updefault}$a$}}}}}
\put(2472,825){\makebox(0,0)[lb]{\smash{{{\SetFigFont{12}{14.4}{\rmdefault}{\mddefault}{\updefault}$b$}}}}}
\end{picture}
\end{center}
\caption{Moore diagram for the lamplighter
automaton~\label{Moorediagram}}
\end{figure}

Automata are usually represented by Moore diagrams. The Moore
diagram for $\mathcal A$ is a directed graph with vertex set $Q$.
The edges are of the form $q\arr{a\mid q(a)} q_a$.
Figure~\ref{Moorediagram} gives the Moore diagram for a certain
two-state automaton studied by Grigorchuk and \.{Z}uk~\cite{GZlamp}.

It is sometimes convenient to define, for $q\in Q$, the \emph{state
function} $\lambda_q:A\to A$ given by \[\lambda_q(a) = q(a) =
\lambda(q,a)\]  If, for each $q\in Q$, the state function
$\lambda_q$ is a permutation, that is belongs to $S_A$, then one can
easily verify that each state $q$ computes a permutation of
$A^*$~\cite{GNS,selfsimilar}.      We call such an automaton
\emph{invertible}.  In particular, if the alphabet of the invertible
automaton is $A_k$ and $q$ is a state, then the function $q$ belongs
to $\mathrm{Aut}(T_k) = S_k\wr \mathrm{Aut}(T_k)$. The wreath
product coordinates of $q$ are:
\begin{equation}\label{wreathcoords}
q = \lambda_q(q_0,\ldots,q_{k-1})
\end{equation}
and so our two uses of the notations $\lambda_q$ and $q_i$ are
consistent. For instance, the automaton from
Figure~\ref{Moorediagram} is described in wreath product coordinates
by $a=(a,b)$, $b=(01)(a,b)$. More generally, if $w\in A_k^*$, then
the section of $q$ at $w$ is exactly the state $q_w$ and in
particular the transformation $q$ is finite state.  One can
show~\cite{GNS,selfsimilar} that the inverse of $q$ is given by the
finite state automaton obtained by switching the two sides of the
labels of the Moore diagram and choosing as the initial state the
state corresponding to $q$.  If $\Ac$ is an invertible automaton,
then $\GG(\Ac)$ denotes the group generated by the states of $\Ac$.
Such groups are called \emph{automaton groups} and constitute the
main examples of finitely generated self-similar
groups~\cite{selfsimilar}. For instance the group generated by the
states of the automaton in Figure~\ref{Moorediagram} is the
lamplighter group $\bigoplus_{\mathbb
Z}\mathbb{Z}/2\mathbb{Z}\rtimes \mathbb Z$ \cite{GNS,GZlamp,reset}.

If $f\in \mathrm{Aut}(T_k)$ is finite state, then it can be computed
by the initial automaton whose state set is $Q=\{f_w\mid w\in A^*\}$
(note: this set is finite by assumption).  The transition and output
functions are given by $\delta(f_w,a) = f_{wa}$ and $\lambda(f_w,a)
= f_w(a)$. The initial state is $f_{\varepsilon}=f$.   We remark
that the composition of finite state transformations is also finite
state~\cite{Eilenberg,GNS,selfsimilar} and so the collection of
invertible finite state maps is a subgroup of $\mathrm{Aut}(T_k)$.

If $H$ is a profinite group, we denote by $[H,H]$ the \emph{closure}
of the commutator subgroup of $H$.  The abelianization $H/[H,H]$ of $H$
shall be denoted $H^{ab}$ and is again a profinite group. Let $(G,A_k)$ be a
transitive permutation group. Then the infinite
permutational wreath product
\begin{equation}\label{infinitewreath}
\wh G = \wr^{\infty} (G,A_k) = (G,A_k)\wr (G,A_k)\wr \cdots
\end{equation}
 is a closed subgroup of
$\mathrm{Aut}(T_k)$. Moreover, it acts spherically transitively on
$T_k$~\cite{Bass}.  The abelianization $\wh{G}^{ab}$ is well known
to be isomorphic to the infinite direct product $G^{ab}\times
G^{ab}\times \cdots$~\cite[Chapter 4, Proposition 4.3]{Bass}. To
describe the map, we think about $\wh{G}^{ab}$ in a different way.
Since $G^{ab}$ is a finite abelian group, it is a finite direct
product of cyclic groups of prime power order in a unique way. Hence
we can view it as the additive group of a finite commutative ring
via this decomposition.  In particular, if $G^{ab}$ is cyclic of
prime order $p$, we view it as the additive group of the field of
$p$ elements. We can then identify $\wh{G}^{ab}$ with the additive
group of the ring of formal power series $G^{ab}\power t$ over
$G^{ab}$ in a single variable $t$.  If $s\in G^{ab}\power t$, we use
the notation $\inner s {t^n}$ to denote the coefficient of $t^n$ in
$s$.  The abelianization map, with this notation, is given
by~\cite{Bass}:
\begin{equation}\label{abelianization}
\inner {g[\wh G,\wh G]} {t^n} = \sum_{|w|=n} \lambda_{g_w}[G,G]
\end{equation}
The importance of the abelianization map is reflected in the
following theorem~\cite[Chapter 4, Propositions (4.6) and
(4.7)]{Bass}.

\begin{theorem}[\cite{Bass}]\label{bigbass}
Let $\wh G = \wr^{\infty}(\mathbb Z/k\mathbb Z,A_k)$.  Then:
\begin{enumerate}
\item  an element $g\in \wh G$ is spherically transitive if and only if its
abelianization $g[\wh G,\wh G]\in \mathbb Z/k\mathbb Z\power t$
satisfies $\inner {g[\wh G,\wh G]} {t^n}\in \mathbb Z/k\mathbb
Z^{\times}$, for all $n\geq 0$; \item two spherically transitive
elements $f,g\in \wh G$ are conjugate if and only if they have the
same image in $\wh G^{ab}=\mathbb Z/k\mathbb Z\power t$.
\end{enumerate}
\end{theorem}

We sketch a proof of the first part of the theorem.   The proof goes
by induction on the levels of the tree and we merely illustrate how
the inductive step works. The key point is that $\langle g\rangle$
acts transitively on $A_k^n$ if and only if it acts transitively on
$A_k^{n-1}$ and, for each word $u\in A_k^{n-1}$, the stabilizer of
$u$ acts transitively on $uA_k$. Now if we assume that $g$ acts as a
$k^{n-1}$-cycle $\sigma$ on $A_k^{n-1}$, then  $g^{k^{n-1}}$
generates the stabilizer of every word in $A_k^{n-1}$. Using the
iterated wreath product decomposition, we can write
$g=\sigma(g_{w_1},\ldots, g_{w_{k^n}})$ where $A_k^{n}
=\{w_1,\ldots,w_{k^n}\}$.  A straightforward calculation then shows
that $g^{k^{n-1}} = (h_1,\ldots, h_{k^n})$ where $h_i
=g_{w_{i-1}}g_{w_{i-2}}\cdots g_{w_1} g_{w_{k^n}}g_{w_{k^n-1}}\cdots
g_{w_i}$.  In particular, $\lambda_{h_i} = \sum_{|w|=n}\lambda_{g_w}
= \inner{g[\wh G,\wh G]} {t^n}$, for all $i$. It follows that
$g^{k^{n-1}}$ acts transitively on $uA_k$ for all $u\in A_k^{n-1}$
if and only if $\inner{g[\wh G,\wh G]} {t^n}\in \mathbb Z/k\mathbb
Z^{\times}$.

Let us return to the setting where $(G,A_k)$ is a transitive
permutation group and let $\wh G$ be as in \eqref{infinitewreath}.
It is easy to see from \eqref{wreathcoords} that if
$\Ac=(Q,A_k,\delta,\lambda)$ is a finite state automaton, then
$\GG(\Ac)\leq \wh G$ if and only if $\lambda_q\in G$ for all $q\in
Q$.

We are now in a position to present our results.  Our first result
provides a partial solution to a problem of Grigorchuk,
Nekrashevych, Sushchanskii and \v Suni\'k from~\cite{GNS}
and~\cite{selfsimilarsurvey}.

\begin{theorem}\label{main1}
Let $g\in \wr^{\infty}(\mathbb Z/k\mathbb Z,A_k)$ be a finite state
transformation, given by a finite state initial automaton. Then it
is decidable whether $f$ is spherically transitive.
\end{theorem}

Our second theorem concerns conjugacy of finite state elements.

\begin{theorem}\label{main2}
Let $f,g\in\wh{G}= \wr^{\infty}(\mathbb Z/k\mathbb Z,A_k)$ be finite
state transformations, given by finite state initial automata. Then
it is decidable whether $f$ and $g$ are conjugate in
$\wh{G}$.\end{theorem}

Theorem~\ref{main2} can be deduced from Theorem~\ref{bigbass} and
our next theorem.

\begin{theorem}\label{main3}
Let $(G,A_k)$ be a transitive permutation group and let $\wh
G=\wr^{\infty}(G,A_k)$.  Let $f,g\in \wh G$ be finite state
transformations, given by finite state initial automata.  Then it is
decidable whether $f$ and $g$ are equal in $\wh G^{ab}$.
\end{theorem}

The key idea for proving these results was inspired by
Sch\"utzenberger's theory of automata and rational power
series~\cite{schutz,schutz2}.  In fact, a biproduct of the proofs
is:

\begin{theorem}\label{main4}
Let $(G,A_k)$ be a transitive permutation group and let $\wh G=
\wr^{\infty}(G,A_k)$.  Let $f\in \wh G$ be a finite state
transformation. Then $f[\wh G,\wh G]\in G^{ab}\power t$ is a
rational power series.
\end{theorem}

\section{Proofs of the theorems}
All the theorems rely on the following simple lemma.
\begin{lemma}\label{matrixthing}
Let $(G,A_k)$ be a transitive permutation group and let $\wh G$ be
as in \eqref{infinitewreath}.  Let $g\in \wh G$ be computed by an
automaton $\Ac$ with state set $\{1,\ldots,n\}$ and initial state
$1$.  Let $A$ be the incidence matrix of $\Ac$ and let $v_A$ be the
vector whose entries are given by $(v_A)_i = \lambda_i[G,G]$,
$i=1,\ldots,n$. Then
\[g[\wh G, \wh G] = \sum_{j=0}^{\infty} (A^jv_A)_1t^j\]
\end{lemma}
\begin{proof}
As $(A^j)_{rs}$ counts the number of paths in $\Ac$ of length $j$
from $r$ to $s$:
\[(A^jv_A)_1 = \sum_{|w|=j}\lambda_{1_w}[G,G] =
\sum_{|w|=j}\lambda_{g_w}[G,G]= \inner {g[\ov G,\ov G]} {t^j}\]
where the last equality follows from \eqref{abelianization}.
\end{proof}

\subsubsection*{Proof of Theorem~\ref{main1}}
 By Theorem~\ref{bigbass}, it follows that $g$
is spherically transitive if and only if each coefficient of $g[\wh
G,\wh G]$ belongs to $\mathbb Z/k\mathbb Z^{\times}$.  By
Lemma~\ref{matrixthing}, we thus want to check whether (keeping the
above notation) $(A^jv_A)_1\in \mathbb Z/k\mathbb Z^{\times}$ for
each $j\geq 0$. Since $\mathbb Z/k\mathbb Z^n$ has $k^n$ elements,
$A^rv_1 =A^sv_1$ for some $0\leq r<s\leq k^n$ and so the above
condition is a finite check.\qed

\subsubsection*{Proof of Theorem~\ref{main3}}
Let $(G,A_k)$ be a transitive permutation group and let $\wh G$ be
as in \eqref{infinitewreath}.  Let $\Ac$ and $\Bc$ be initial
automata computing $f$ and $g$, respectively.  Say that $\Ac$ has
$m$ states and $\Bc$ has $n$ states. Let $A$ and $B$ be the
respective incidence matrices of $\Ac$ and $\Bc$.  Let $v_A$ and
$v_B$ be the associated vectors, as per Lemma~\ref{matrixthing}.
Consider the matrix $M = \begin{pmatrix} A
  & 0\\ 0 &
  B\end{pmatrix}$. Let $\{e_1,\ldots, e_{m+n}\}$ be
the standard basis
of row vectors for $(G^{ab})^{m+n}$ and set $v=\begin{pmatrix}
  v_A\\ v_B\end{pmatrix}$. Then, applying Lemma~\ref{matrixthing},
  we have for $j\geq 0$: \[(e_1-e_{m+1})(M^jv) =
(A^jv_A)_1 - (B^jv_B)_1 = \inner {f[\wh G,\wh G]} {t^j}-\inner
{g[\wh G,\wh G]} {t^j}\] Hence
  $f[\wh G,\wh G]=g[\wh G,\wh G]$ if and only if
  $(e_1-e_{m+1})(M^jv)=0$ for all $j\geq 0$.  But again, $M^rv=M^sv$ some $0\leq r<s\leq k^{m+n}$,
   so we can check this.

   If $G^{ab}$ is a finite
  field, then we can do better.  Indeed, since the vectors
  $v,Mv,\ldots, M^{m+n}v$ in $\mathbb Z/k\mathbb Z^{m+n}$ must be linearly dependent, it follows that for some
   $0\leq i\leq m+n$, $M^{i}v = c_0v+ c_1Mv\cdots+ c_{i-1}M^{i-1}v$.
   Such a recursion implies that $M^jv$ is a linear combination of
   $v,Mv,\ldots, M^{n+m-1}v$ for all $j\geq n+m$.  Hence $(e_1-e_{m+1})(M^jv)
   =0$ for all $j\geq 0$ if and only if $(e_1-e_{m+1})(M^jv)
   =0$ for $0\leq j\leq m+n-1$.
\qed

\begin{remark}
The proof of Theorem~\ref{main3} allows for an alternative algorithm
for testing spherical transitivity for $\mathrm{Aut}(T_2)$. By
Theorem~\ref{bigbass}, $g\in \mathrm{Aut}(T_2)$ is spherically
transitive if and only if $g[\mathrm{Aut}(T_2),\mathrm{Aut}(T_2)] =
\sum_{n=0}^{\infty}t^n$, and all spherically transitive elements are
conjugate.  The so-called \emph{odometer} $a=(01)(1,a)$ is one such
spherically transitive element and it has two distinct sections,
that is, it can be computed by a two-state automaton.  It follows
from the proof of Theorem~\ref{main3} that if $g\in
\mathrm{Aut}(T_2)$ is computed by an $n$-state initial automaton
with incidence matrix $A$, then one needs only to verify
$(A^jv_A)_1\neq 0$ for $0\leq j\leq n+1$.
\end{remark}

\subsubsection*{Proof of Theorem~\ref{main4}}
From Lemma~\ref{matrixthing} that we have $g[\wh G,\wh G] =
((I-At)^{-1}v_A)_1$. Since \[(I-At)\inv =
\frac{1}{\mathrm{det}(I-At)}\mathrm{Adj}(I-At)\] and the entries of
the adjoint $\mathrm{Adj}(I-At)$ are polynomials in $t$, while
\mbox{$\mathrm{det}(I-At)$} is a polynomial in $t$, it follows that
the entries of $(I-At)\inv$ are rational power series in $t$. Since
$((I-At)^{-1}v_A)_1$ is a linear combination of entries of
$(I-At)^{-1}$, it follows that  $g[\wh G,\wh G]$ is a rational pwer
series.\qed

\section*{Acknowledgments}
We would like to thank Zoran \v Suni\'k for some helpful comments
and his careful reading of an earlier draft of this paper.

\bibliographystyle{amsplain}

\end{document}